\documentclass[MSNbibl,number,citesort,seceqn,dvips]{arxbj}
\usepackage{algorithm}
\usepackage{graphicx}

\aid{0}
\volume{18}
\issue{4}
\pubyear{2012}
\firstpage{1320}
\lastpage{1340}
\doi{10.3150/11-BEJ380}

\makeatletter
\newtheorem{theo}{Theorem}[section]
\newtheorem{lemma}{Lemma}[section]
\newcommand{\fraca}[2]{{#1}\big/{#2}}
\makeatother

\begin{document}
\begin{frontmatter}

\title{On a characterization of ordered pivotal sampling}
\runtitle{On a characterization of ordered pivotal sampling}

\begin{aug}
%%%% inicialai - be tarpu
\author{\fnms{Guillaume} \snm{Chauvet}\corref{}\ead[label=e1]{chauvet@ensai.fr}}% \and
\runauthor{G. Chauvet} %% auto
\address{Crest (Ensai), Campus de Ker Lann, 35170 Bruz, France. \printead{e1}}
\end{aug}

% HISTORY:
\received{\smonth{11} \syear{2010}}
\revised{\smonth{4} \syear{2011}}

% ABSTRACT
%
\begin{abstract}
When auxiliary information is available at the design stage, samples
may be selected by means of balanced sampling. Deville and Till\'e
proposed in 2004 a general algorithm to perform balanced sampling,
named the cube method. In this paper, we are interested in a particular
case of the cube method named pivotal sampling, and first described by
Deville and Till\'e in 1998. We show that this sampling algorithm, when
applied to units ranked in a fixed order, is equivalent to Deville's
systematic sampling, in the sense that both algorithms lead to the same
sampling design. This characterization enables the computation of the
second-order inclusion probabilities for pivotal sampling. We show that
the pivotal sampling enables to take account of an appropriate ordering
of the units to achieve a variance reduction, while limiting the loss
of efficiency if the ordering is not appropriate.
\end{abstract}

% KEYWORDS
%

\end{frontmatter}

%s1 ###
%se1 #&#
\section{Introduction}

When auxiliary information is available at the design stage, samples
may be selected by means of balanced sampling. The variance of the
Horvitz--Thompson (HT) estimator is then reduced, since it is
approximately given by that of the residuals of the variable of
interest on the balancing variables.
Deville and Till\'e~\cite{devtil04} proposed a general algorithm for balanced sampling,
named the \textit{cube method}. This sampling algorithm enables the
selection of balanced samples with any number of balancing variables,
and any prescribed set of inclusion probabilities.

In order to measure the gain in efficiency provided by the cube method,
Deville and Till\'e~\cite{devtil05} proposed several variance approximations. They
suppose that the sampling design is exactly balanced, and performed
with maximum entropy among sampling designs balanced on the same
balancing variables, with the same inclusion probabilities. Then, under
an additional assumption of asymptotic normality of the multivariate
HT-estimator under Poisson sampling, the variance approximations are
derived. The assumption of exact balancing may be closely respected, if
the number of balancing variables remains small with regard to the
sample size; otherwise, the balancing error must be taken into account
in variance estimation, see Breidt and Chauvet~\cite{brecha11}. The second assumption is
related to the entropy of the sampling design: the variance
approximations proposed by Deville and Till\'e~\cite{devtil05} are unlikely to hold if
this assumption is not satisfied.

A practical way to increase the entropy of a sampling design is to sort
the population randomly before sampling. However, this preliminary
randomization step is not systematically included in the sampling
process. This is a common practice to sort the population with respect
to some auxiliary variable before the sampling, so as to benefit from a
stratification effect. In France, Census surveys are conducted
annually; the detailed methodology is described in Godinot~\cite{god05}. Each
large municipality (10\,000 inhabitants or more in 1999) is the subject
of an independent sampling design and is stratified according to the
type of address (large addresses, new addresses, or other addresses).
In each stratum, the addresses are divided into 5 rotation groups. Each
year, all the addresses within one rotation group (for the strata of
large addresses and new addresses) or within a sub-sample (for the
stratum of other addresses) are surveyed. In the stratum of other
addresses, the sub-sample is obtained by first, sorting the addresses
with respect to the descending number of dwellings, and then, applying
the cube method. In such cases, the conditions for the variance
approximations proposed by Deville and Till\'e~\cite{devtil05} to hold are clearly not respected.

We are interested in a particular case of the cube method, called \textit{pivotal sampling} (Deville and Till\'e~\cite{devtil98}), obtained when the only balancing
condition is given by the variable of inclusion probabilities. That is,
the cube method with the sole fixed-size constraint amounts to pivotal
sampling. This algorithm is an exact sampling procedure, which respects
a prescribed set of inclusion probabilities, is strictly without
replacement and leads to fixed-size designs. In this paper, we show
that the pivotal sampling algorithm, when applied to units ranked in a
fixed order, is equivalent to an algorithm proposed in Deville~\cite{dev88},
and known in the literature as \textit{Deville's systematic sampling}
(Till\'e~\cite{til06}). The two algorithms are equivalent, in the sense that
both lead to the same sampling design. In particular, the computation
of the second-order inclusion probabilities developed in Deville~\cite{dev88}
may be readily applied to pivotal sampling. This provides an answer to
a problem raised by  Bondesson and Grafstr\"om~\cite{bongra10}, page 7. Deville's systematic
sampling has similarities with Markov chain designs introduced by  Breidt~\cite{bre95}. It has found uses in the context of longitudinal surveys, see
Nedyalkova, Qualit\'e, and Till\'e~\cite{nedquatil09}.

The paper is organized as follows. In Section~\ref{secnot}, the
notation is defined. Ordered pivotal sampling and Deville's systematic
sampling are presented in Sections~\ref{secops} and~\ref{secdss},
respectively, and some useful results are derived. The second-order
inclusion probabilities for ordered pivotal sampling are given in
Section~\ref{secsoip}. Some results which illustrate the practical
interest of ordered pivotal sampling are presented in Section~\ref{intops}.

%s2 ###
%se2 #&#
\section{Notation} \label{secnot}

Consider a finite population $U$ consisting of $N$ sampling units that
may be represented by integers $k=1,\ldots,N$. We assume that the order
of the units in the population is fixed prior\vadjust{\goodbreak} to sampling, and may be
confounded with the natural order of their indexes. A sample $s$,
defined as a subset of $U$, is selected with inclusion probabilities
$\pi= (\pi_1,\ldots,\pi_N )'.$ We assume without loss of
generality that $0 < \pi_k <1$ for any unit $k$ in $U$, with $n=\sum_{k
\in U} \pi_k$ the sample size. Let $\pi_{kl}$ denote the probability
that units $k$ and $l$ are selected jointly in the sample.

\begin{figure}

\includegraphics{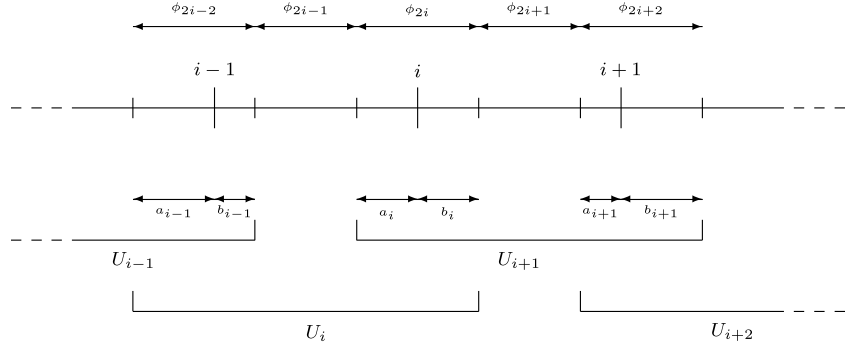}

\caption{Inclusion probabilities and cross-border units in
microstratum $U_i$, for population $U$.} \label{nota}
\end{figure}

We define $V_k=\sum_{l=1}^k \pi_l$ for any unit $k \in U$, with
$V_0=0$. A unit $k$ is said to be \textit{cross-border} if $V_{k-1} \le i$
and $V_{k} > i$ for some nonnegative integer $i$. The cross-border
units are denoted as $k_i, i=1,\ldots,n-1,$ and we note
$a_i=i-V_{k_i-1}$ and $b_i=V_{k_i}-i$. The microstratum $U_i, i=1,\ldots
,n,$ is defined as
%
%e2.1 ###
%e2.1 #&#
\begin{equation} \label{Ui}
U_i=\{k \in U; k_{i-1} \leq k \leq k_i\},
\end{equation}
with $k_0=0$ and $k_n=N+1$. To fix ideas, useful quantities for
population $U$ are presented in Figure~\ref{nota}.

%f1 ###
%fi1 #&#

The microstrata are generally overlapping, since one cross-border unit
may belong to two adjacent microstrata: the cross-border unit $k_i$
belongs both to the microstratum $U_i$ (with an associated probability
$a_i$) and to the microstratum $U_{i+1}$ (with an associated
probability $b_i$). In the particular case when $V_{k_i}=i$, we have
$b_i=0$. To avoid the introduction of specific notations for such
cases, we consider in Sections~\ref{secnot}--\ref{secsoip} that, in
this situation, the cross-border unit $k_i$ belongs to the microstratum
$U_{i+1}$ as a ``phantom unit,'' that is, with an associated probability
equal to 0. In Section~\ref{intops}, we simply consider that the
cross-border unit $k_i$ belongs to the microstratum $U_i$ only in such
situations.

\begin{figure}

\includegraphics{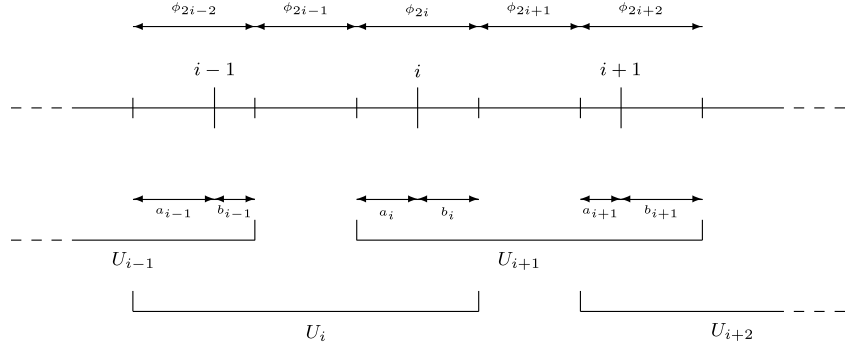}

\caption{Inclusion probabilities and cross-border units in microstrata
$U_{i}$ and $U_{i+1}$
for population $U_c$.} \label{notaclust}
\end{figure}

The $N$ sampling units are grouped to obtain a population
$U_c=\{u_1,\ldots,u_{2n-1}\}$ of clusters. There are the clusters of
cross-border units ($n-1$ singletons), denoted as $u_{2i}$ with
associated probability $\phi_{2i}=\pi_{k_i}$ for $i=1,\ldots,n-1$.
There are the $n$ clusters of units that are not cross-borders and
that are between two consecutive integers, denoted as $u_{2i-1}$
with associated probability $\phi_{2i-1}=V_{k_i-1}-V_{k_{i-1}}$, for
$i=1,\ldots,n$. We note
$\psi= (\phi_1,\ldots,\phi_{2n-1} )'$. To fix ideas, useful
quantities for population $U_c$ are presented in Figure
\ref{notaclust}. If (at least) one of the cross-border units in
$U_i$ has a large inclusion probability, there may not exist any non-cross-border unit between integers $i-1$ and $i$, so that the
cluster $u_{2i-1}$ is empty. To avoid the need for specific
notations for such cases, we may view this situation as a particular
case of our framework by allowing a cluster $u_{2i-1}$ to be a
``phantom cluster,'' that is, an empty cluster with associated
probability $\phi_{2i-1}$ equal to~$0$. For example, suppose that
$N=8$, $n=4$ and $\pi=(0.2,0.5,0.3,0.4,0.9,0.8,0.5,0.4)'$. We obtain
the $4$ microstrata $U_1=\{1,2,3\}$, $U_2=\{3,4,5\}$, $U_3=\{5,6\}$
and $U_4=\{6,7,8\}$. In particular, we have $a_1=0.3=\pi_3$ and
$b_1=0$, so that the cross-border unit $3$ is a phantom unit for the
microstratum $U_2$. Also, we obtain $7$ clusters (see Table
\ref{exclust}): the cluster $u_5$ is empty, with an associated
probability equal to zero.

%s3 ###
%se3 #&#
\section{Ordered pivotal sampling} \label{secops}

A general algorithm for pivotal sampling is described in  Deville and Till\'e~\cite{devtil98}. In the version presented in Algorithm~\ref{pivmet}, the
order of the sampling units is explicitly taken into account. We call
it \textit{ordered pivotal sampling} to avoid confusion. At each step, one
or two coordinates of $\pi(t)$ are randomly rounded to $0$ or $1$, and
remain there forever. In at most $N$ steps, the final sample is obtained.

%
%f2 ###
%fi2 #&#

%t1 ###
%ta1 #&#
\begin{table*}[b]
\tabcolsep=0pt
\tablewidth=295pt
\caption{Clusters and associated probabilities for a population of size~$8$} \label{exclust}
\begin{tabular*}{295pt}{@{\extracolsep{\fill}}llllllll@{}}
\hline
$i$ & 1 & 2 & 3 & 4 & 5 & 6 & 7 \\
\hline
$u_i$ & $\{1,2\}$ & $\{3\}$ & $\{4\}$ & $\{5\}$ & $\emptyset$ & $\{6\}
$ & $\{7,8\}$ \\
$\phi_i$ & 0.7 & 0.3 & 0.4 & 0.9 & 0 & 0.8 & 0.9 \\
\hline
\end{tabular*}
\end{table*}

Roughly speaking, the algorithm may be summarized as follows. At the
beginning, in microstratum $U_1$ ($i=1$), the two first units $1$ and
$2$ fight, the loser is definitely eliminated while the survivor
(denoted as $J_0$) gets the sum of their probabilities and then faces
the following unit. The fights go on until the accumulated probability
exceeds $1$, which occurs at time $t=k_1$ when the survivor $J_0$ faces
the cross-border unit $k_1$. One of the two remaining units, denoted as
$W_1$, wins and is then definitely selected in the sample while the
other one, denoted as $J_1$, jumps to the microstratum $U_2$.

More generally, in microstratum $U_i$, the first unit $k_{i-1}$ is
replaced with the unit $J_{i-1}$ which jumps from the microstratum
$U_{i-1}$. The units $J_{i-1}$ and $k_{i-1}+1$ fight, the survivor gets
the sum of their probabilities and then faces the next unit. The fights
go on until the survivor $J_{i-1}$ faces the cross-border unit $k_i$.
One of the two remaining units ($W_i$) wins and is then definitely
selected in the sample while the other one ($J_i$) jumps to the
following microstratum. Lemma~\ref{twostage} states that Algorithm \ref
{pivmet} may alternatively be seen as a two-stage procedure. The proof
follows from definition, and is thus omitted.

\begin{algorithm}[t]
\begin{enumerate}[3.]
\item We initialize with $i=1$, $J_0=1$ and $\pi(0)=\pi$.
\item For $t=2,\ldots,N,$ do:
\begin{enumerate}[(b)]
\item[(a)] If $m \in U \setminus\{J_{i-1},t\}$, then $\pi_m(t)=\pi_m(t-1)$.
\item[(b)] If $\pi_{J_{i-1}}(t-1)+\pi_t(t-1)<1$, let
$\lambda_1(t)=\frac{\pi_{J_{i-1}}(t-1)}{\pi_{J_{i-1}}(t-1)+\pi_t(t-1)}.$
Then
\begin{enumerate}[ii.]
\item with probability $\lambda_1(t)$, let
\[
 [\pi_{J_{i-1}}(t),\pi_t(t)  ]= [\pi_{J_{i-1}}(t-1)+\pi_t(t-1),0
 ];
\]
\item with probability $1-\lambda_1(t)$, let $J_{i-1}=t$
and
\[
 [\pi_{J_{i-1}}(t),\pi_t(t)  ]= [0,\pi_{J_{i-1}}(t-1)+\pi_t(t-1)
 ].
\]
\end{enumerate}
\item[(c)] If $\pi_{J_{i-1}}(t-1)+\pi_t(t-1) \geq1$, let $\lambda_1(t)=\frac
{1-\pi_t(t-1)}{2-\pi_{J_{i-1}}(t-1)-\pi_t(t-1)}.$
Then
\begin{enumerate}[iii.]
\item with probability $\lambda_1(t)$, let $W_i=J_{i-1}$, let $J_{i}=t$ and
\[
 [\pi_{J_{i-1}}(t),\pi_t(t)  ]= [1,\pi_{J_{i-1}}(t-1)+\pi_t(t-1)-1
 ];
\]
\item with probability $1-\lambda_1(t)$, let
$W_i=t$, let $J_{i}=J_{i-1}$ and
\[
 [\pi_{J_{i-1}}(t),\pi_t(t)  ]= [\pi_{J_{i-1}}(t-1)+\pi_t(t-1)-1,1
 ];
\]
\item let $i=i+1$.
\end{enumerate}
\end{enumerate}
\item The sample is given by $\{W_1,\ldots,W_n\}$.
\end{enumerate}
\caption{Ordered Pivotal Sampling with parameter $\boldsymbol{\pi}$}
\label{pivmet}
\end{algorithm}

%le3.1 #&#
\begin{lemma} \label{twostage}
Ordered pivotal sampling with parameter $\pi$ may be obtained by
two-stage sampling, where a sample
$s_c$ of $n$ clusters is first selected in $U_c$ by means of ordered
pivotal sampling with parameter $\psi$,
and one unit $k$ is then selected in each $u_j \in s_c$ with a
probability proportional to $\pi_k$.
%a sample of size $1$ is then selected in each $u_j \in s_c$ with
%probabilities proportional to $\pi_k$.
\end{lemma}

We assume that a sample $S_{op}$ is selected in $U_c$ by means of
ordered pivotal sampling with parameter $\psi$, and we let $X_1<\cdots <X_n$ denote the units selected in the sample, ranked in ascending
order. Lemma~\ref{equwinjum} states useful relations between on the
one hand, the sampled units $X_i$, and on the other hand, the winners
$W_i$ and jumpers $J_i$. Lemma~\ref{probu2i-1} gives the probabilities
for the different outcomes in the case of a non-cross-border unit~$u_{2i-1}$.

%le3.2 #&#
\begin{lemma} \label{equwinjum}
In case of ordered pivotal sampling with parameter $\psi$, we have
%
%e3.3 ###
%e3.2 ###
%e3.1 ###
%e3.1 #&#
%e3.2 #&#
%e3.3 #&#
\begin{eqnarray} \label{equwj1}
\{X_i=u_{2i-2}\} & \Rightarrow& \bigl\{J_{i-1} \in\{ X_1,\ldots,X_i \} \bigr\}
, \\\label{equwj2}
\{X_i=u_{2i-1}\} & \Rightarrow& \{W_i=u_{2i-1}\} \cup\{J_i=u_{2i-1}\}
,  \\
\label{equwj3}
\{X_i=u_{2i}\} & \Rightarrow& \bigl\{J_i \notin\{ X_1,\ldots,X_i \} \bigr\}.
\end{eqnarray}
\end{lemma}

\begin{pf}%[\ref{equwinjum}]
Assume that $X_i=u_{2i-2}$. This implies that $i$ units exactly are
selected in the $i-1$ first microstrata $U_1,\ldots,U_{i-1}$. On the
other hand, if $J_{i-1} \notin\{ X_1,\ldots,X_i \}$ the unit
$J_{i-1}$ is not selected in the sample so that at most $i-1$ units
are selected in $U_1,\ldots,U_{i-1}$. This proves (\ref{equwj1}),
and by a similar argument we obtain (\ref{equwj3}). It is easily
seen that (\ref{equwj2}) holds, since the selection of $u_{2i-1}$
implies that this unit is either the winner $W_i$ or the jumper
$J_i$ in the microstratum~$U_i$.
\end{pf}

%le3.3 #&#
\begin{lemma} \label{probu2i-1}
In case of ordered pivotal sampling with parameter $\psi$, we have
%
%e3.6 ###
%e3.5 ###
%e3.4 ###
%e3.4 #&#
%e3.5 #&#
%e3.6 #&#
\begin{eqnarray}\label{probu2i-1w}
pr  ( W_i=u_{2i-1} ) & = & \frac
{(1-a_i-b_{i-1})(1-a_i-b_i)}{(1-a_i)(1-b_i)},  \\\label{probu2i-1j}
pr  ( J_i=u_{2i-1} ) & = & \frac
{a_i(1-a_i-b_{i-1})}{(1-a_i)(1-b_i)},  \\\label{probu2i-1x}
pr  ( X_i=u_{2i-1} ) & = & 1-a_i-b_{i-1}.
\end{eqnarray}
\end{lemma}

\begin{pf}%[\ref{equwinjum}]
The event
\[
\{ W_i=u_{2i-1} \}
\]
may be alternatively interpreted as follows: in the fight
between $J_{i-1}$ and $u_{2i-1}$, the unit $u_{2i-1}$ survives;
then in the next fight, the unit $u_{2i-1}$ is the selected unit
$W_i$, while the unit $u_{2i}$ is the jumping unit $J_i$.
Consequently, we have:
\begin{eqnarray*}
pr  ( W_i=u_{2i-1} ) & = &
\frac{1-b_{i-1}-a_i}{1-a_i} \times
\frac{1-a_i-b_i}{1-b_i},
\end{eqnarray*}
which gives (\ref{probu2i-1w}). Similarly, we obtain
\begin{eqnarray*}
pr  ( J_i=u_{2i-1} ) & = &
\frac{1-b_{i-1}-a_i}{1-a_i} \times
\frac{a_i}{1-b_i},
\end{eqnarray*}
which gives (\ref{probu2i-1j}). We now consider equation
(\ref{probu2i-1x}). Since
\[
\{X_i=u_{2i-1}\} \Rightarrow\{u_{2i-1} \in S_{op}\}
\]
and
\[
pr  ( u_{2i-1} \in S_{op}  )=1-a_i-b_{i-1},
\]
it suffices to show that
%
%e3.7 ###
%e3.7 #&#
\begin{equation} \label{calvin}
\{u_{2i-1} \in S_{op}\} \Rightarrow\{X_i=u_{2i-1}\}.
\end{equation}
Since $\{u_{2i-1} \in S_{op}\}$ implies that $u_{2i-1}$ survives in its
duel against
$J_{i-1}$, this in turn implies that $J_{i-1} \notin\{ X_1,\ldots,X_i
\}$. In other words, $\{u_{2i-1} \in S_{op}\}$ implies that exactly
$i-1$ units smaller than $u_{2i-1}$ were selected, which proves
(\ref{calvin}).
\end{pf}

Finally, let $U_{c,i}=\{u_{2i-2},\ldots,u_{2n-1}\}$, $\psi_i=
(b_{i-1},\phi_{2i-1},\ldots,\phi_j,\ldots,\phi_{2n-1} )'$,
and $S_{op,i}$ be a random sample selected in $U_{c,i}$ by means of
ordered pivotal sampling with parameter $\psi_i$. Lemma~\ref{incprobopsi}
establishes some relations for conditional inclusion probabilities in
$S_{op,i}$ of the first units
in $U_{c,i}$.

%le3.4 #&#
\begin{lemma} \label{incprobopsi}\vspace*{-1.75pt}
%
%e3.10 ###
%e3.9 ###
%e3.8 ###
%e3.8 #&#
%e3.9 #&#
%e3.10 #&#
\begin{eqnarray} \label{opsi1}
& & pr ( u_{2i} \in S_{op,i}, u_{2i-1} \notin S_{op,i} |
u_{2i-2} \in S_{op,i}  )\nonumber
\\[-8pt]
\\[-8pt]
&& \quad  =  \frac{b_i}{1-a_i}, \nonumber\\\label
{opsi2}
& & pr ( u_{2i+1} \in S_{op,i},u_{2i} \notin S_{op,i}, u_{2i-1}
\notin S_{op,i} | u_{2i-2} \in S_{op,i}  )\nonumber
\\[-8pt]
\\[-8pt]
&& \quad  =  \frac{(1-a_i-b_i)(1-b_i-a_{i+1})}{(1-a_i)(1-b_i)},  \nonumber\\\label{opsi3}
& & pr ( u_{2i+2} \in S_{op,i},u_{2i+1} \notin S_{op,i},u_{2i}
\notin S_{op,i}, u_{2i-1} \notin S_{op,i} | u_{2i-2} \in S_{op,i}
 ) \nonumber
 \\[-8pt]
 \\[-8pt]
&& \quad  =  \frac{(1-a_i-b_i)a_{i+1}}{(1-a_i)(1-b_i)}.
\nonumber
\end{eqnarray}
\end{lemma}

\begin{pf}
To fix ideas, the first units in population $U_{c,i}$ and
related quantities are presented in Figure~\ref{nota2}.

%f3 ###
%fi3 #&#
\begin{figure}

\includegraphics{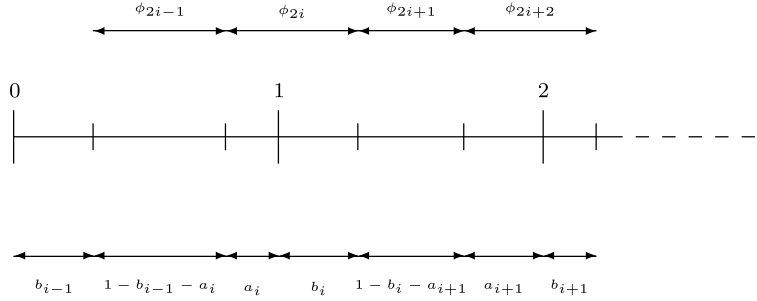}

\caption{Inclusion probabilities and cross-border units in the two
first microstrata of population $U_{c,i}$.} \label{nota2}
\end{figure}

We first consider equation (\ref{opsi1}). Since $b_{i-1}$ is
the first-order inclusion probability of unit $u_{2i-2}$ in
sample $S_{op,i}$, we have\vspace*{-1pt}
%
%e3.11 ###
%e3.11 #&#
\begin{eqnarray} \label{denopsi}
pr ( u_{2i-2} \in S_{op,i}  ) & = & b_{i-1}.
\end{eqnarray}
On the other hand, the event\vspace*{-1pt}
\[
\{u_{2i} \in S_{op,i}, u_{2i-1} \notin S_{op,i}, u_{2i-2} \in
S_{op,i}\}
\]
may be alternatively interpreted as follows: in the first fight,
the unit $u_{2i-2}$ survives against the unit $u_{2i-1}$; in the second
fight, any of the two units $u_{2i-2}$ or $u_{2i}$ is the
selected unit $W_1$, while the other is the jumping unit $J_1$;
then, the jumping unit $J_1$ is selected during one of the following
fights. Consequently, we have:\vspace*{-1pt}
%
%e3.12 ###
%e3.12 #&#
\begin{eqnarray}
\label{num1opsi}
& & pr  ( u_{2i} \in S_{op,i}, u_{2i-1} \notin S_{op,i}, u_{2i-2}
\in
S_{op,i}  ) \nonumber
\\[-8pt]
\\[-8pt]
&& \quad  =  \frac{b_{i-1}}{1-a_i} \times1 \times b_i,
\nonumber
\end{eqnarray}
and equation (\ref{opsi1}) follows from (\ref{denopsi}) and
(\ref{num1opsi}). We now consider equation (\ref{opsi2}). The event
\[
\{u_{2i+1} \in S_{op,i}, u_{2i} \notin S_{op,i}, u_{2i-1} \notin
S_{op,i}, u_{2i-2} \in S_{op,i}\}
\]
may be interpreted as follows: in the first fight,
the unit $u_{2i-2}$ survives against the unit $u_{2i-1}$; in the second
fight, $u_{2i-2}$ is the selected unit $W_1$, while $u_{2i}$ is the
jumping unit~$J_1$; in the third fight, the unit $u_{2i+1}$ survives against the
unit $u_{2i}$;
then, the unit $u_{2i+1}$ is selected during one of the
following fights. Consequently, we have:
%
%e3.13 ###
%e3.13 #&#
\begin{eqnarray}\label{num2opsi}
& & pr  ( u_{2i+1} \in S_{op,i}, u_{2i} \notin S_{op,i}, u_{2i-1}
\notin S_{op,i}, u_{2i-2} \in S_{op,i}  ) \nonumber\\
&& \quad  =  \frac{b_{i-1}}{1-a_i} \times\frac{1-a_i-b_i}{1-b_i} \times\frac
{1-b_i-a_{i+1}}{1-a_{i+1}} \times
(1-a_{i+1}) \\
&& \quad  =
\frac{b_{i-1}(1-a_i-b_i)(1-b_i-a_{i+1})}{(1-a_i)(1-b_i)}, \nonumber
\end{eqnarray}
which, together with (\ref{denopsi}), leads to (\ref{opsi2}).
Finally, we consider equation (\ref{opsi3}). The event
\[
\{u_{2i+2} \in S_{op,i}, u_{2i+1} \notin S_{op,i}, u_{2i} \notin
S_{op,i}, u_{2i-1} \notin S_{op,i}, u_{2i-2} \in S_{op,i}\}
\]
may be interpreted as follows: in the first fight,
the unit $u_{2i-2}$ survives against the unit $u_{2i-1}$; in the second
fight, $u_{2i-2}$ is the selected unit $W_1$, while $u_{2i}$ is the
jumping unit~$J_1$; in the third fight, any of the two units $J_i=u_{2i}$ or
$u_{2i+1}$ survives; in the fourth fight, $u_{2i+2}$ is the selected
unit $W_2$, while the other unit is the
jumper $J_2$; then, the unit $J_2$ is not selected during one of the
following fights. Consequently, we have:
%
%e3.14 ###
%e3.14 #&#
\begin{eqnarray}\label{num3opsi}
& & pr  ( u_{2i+2} \in S_{op,i},u_{2i+1} \notin S_{op,i}, u_{2i}
\notin S_{op,i},
u_{2i-1} \notin S_{op,i}, u_{2i-2} \in S_{op,i}  ) \nonumber\\
&& \quad  =  \frac{b_{i-1}}{1-a_i} \times\frac{1-a_i-b_i}{1-b_i} \times1
\times\frac{a_{i+1}}{1-b_{i+1}} \times
(1-b_{i+1})  \\
&& \quad  =
\frac{b_{i-1}(1-a_i-b_i)a_{i+1}}{(1-a_i)(1-b_i)}, \nonumber
\end{eqnarray}
which gives (\ref{opsi3}).
\end{pf}

\begin{algorithm}[b]
\centerline{At step 1:}
\begin{enumerate}[2.]
\item A distributed Uniform$(0,1)$ random variable $w_1$ is generated.
\item The unit $k$ is selected if $V_{k-1} \leq w_1 < V_k$.
\end{enumerate}
\centerline{At step $i$:}
\begin{enumerate}[2.]
\item A random variable $w_i$ is generated:
\begin{enumerate}[(b)]
\item[(a)] if unit $k_{i-1}$ was selected at step $i-1$, then $w_i$ is
generated according to a distributed Uniform$(b_{i-1},1)$ random
variable,
\item[(b)] otherwise, $w_i$ is generated:
\begin{itemize}
\item according to a distributed Uniform$(0,b_{i-1})$ random variable
with probability $a_{i-1} b_{i-1}  \{(1-a_{i-1})(1-b_{i-1}) \}^{-1}$,
\item according to a distributed Uniform$(0,1)$ random variable with
probability $1-a_{i-1} b_{i-1}  \{(1-a_{i-1})(1-b_{i-1}) \}^{-1}$.
\end{itemize}
\end{enumerate}
\item The unit $k$ is selected if $V_{k-1} \leq w_i+(i-1) < V_{k}$.
\end{enumerate}
\caption{Deville's systematic sampling with parameter $\boldsymbol{\pi
}$} \label{sysdev}
\end{algorithm}

%s4 ###
%se4 #&#
\section{Deville's systematic sampling} \label{secdss}

The sampling algorithm known in the literature as Deville's systematic
sampling (Deville~\cite{dev88}; Till\'e~\cite{til06}) is presented in Algorithm~\ref{sysdev}.
This algorithm proceeds in $n$
sub-samplings of size $1$ in the microstrata $U_1,\ldots,U_n$, and the
random variables $w_i$ which indicate the sampled units are generated
so that a cross-border unit $k_{i-1}$ may not be selected twice in the
sample: at step $i$, one unit denoted as $Y_i$ is drawn in $U_i$ if
$k_{i-1}$ was not selected at step $i-1$, and in $U_i \setminus\{
k_{i-1}\}$ otherwise. This sampling algorithm may be particularly
useful in the context of business surveys, when a fine stratification
is used leading to small and possibly non-integer sample size inside
(micro)strata. Deville's systematic sampling directly handles the
rounding problem, since any unit for which the sampling outcome is
still undecided is moved to the next stratum, where the final sampling
decision is then obtained. Lemma~\ref{twostage2} follows from the
definition of Algorithm~\ref{sysdev}.

%le4.1 #&#
\begin{lemma} \label{twostage2}
Deville's systematic sampling with parameter $\pi$ may be obtained by
two-stage sampling, where a sample
$s_c$ of $n$ clusters is first selected in $U_c$ by means of Deville's
systematic sampling with parameter $\psi$,
and one unit $k$ is then selected in each $u_j \in s_c$ with a
probability proportional to $\pi_k$.
\end{lemma}

Assume that a sample is selected in $U_c$ by means of Deville's
systematic sampling with parameter $\psi$. The random variable
$Y_{i+1}$ which gives the result of the sampling in the microstratum
$U_{i+1}$ only depends on the outcome of step $i$, so that
%
%e4.1 ###
%e4.1 #&#
\begin{equation} \label{markprop}
pr  ( Y_{i+1} = u_j | Y_1,\ldots,Y_i  ) = pr  (Y_{i+1}=u_j | Y_i  ).
\end{equation}
The different cases for the transition probabilities in (\ref
{markprop}) easily follow from the definition of Algorithm \ref
{sysdev}, and are given below:
%e4.2 ###
%e4.2 #&#
%e4.3 #&#
%e4.4 #&#
\begin{eqnarray}
\label{case1}
& & pr  ( Y_{i+1}=u_j | Y_1,\ldots,Y_{i-1},Y_i=u_{2i-2}  )
\nonumber
\\[-8pt]
\\[-8pt]
&& \quad  =\everymath{\displaystyle }
\cases{
\frac{b_i}{1-a_i}, & $j=2i$,\vspace*{2pt}\cr
\frac{(1-b_i-a_{i+1})(1-a_i-b_i)}{(1-a_i)(1-b_i)}, & $j=2i+1$,\vspace*{2pt}\cr
\frac{a_{i+1}(1-a_i-b_i)}{(1-a_i)(1-b_i)}, & $j=2i+2$,
}
\nonumber
\\
\label{case2}
& & pr  ( Y_{i+1}=u_j | Y_1,\ldots,Y_{i-1},Y_i=u_{2i-1}  )
\nonumber
\\[-8pt]
\\[-8pt]
&& \quad  =\everymath{\displaystyle }
\cases{
\frac{b_i}{1-a_i}, & $j=2i$,\vspace*{2pt}\cr
\frac{(1-b_i-a_{i+1})(1-a_i-b_i)}{(1-a_i)(1-b_i)}, & $j=2i+1$,\vspace*{2pt}\cr
\frac{a_{i+1}(1-a_i-b_i)}{(1-a_i)(1-b_i)}, & $j=2i+2$,
}
\nonumber
\\
\label{case3}
& & pr  ( Y_{i+1}=u_j | Y_1,\ldots,Y_{i-1},Y_i=u_{2i}  )
\nonumber
\\[-8pt]
\\[-8pt]
&& \quad  =\everymath{\displaystyle }
\cases{
\frac{(1-b_i-a_{i+1})}{(1-b_i)}, & $j=2i+1$,\vspace*{2pt}\cr
\frac{a_{i+1}}{(1-b_i)}, & $j=2i+2$.
}
\nonumber
\end{eqnarray}

%s5 ###
%se5 #&#
\section{Second-order inclusion probabilities} \label{secsoip}

We can now formulate our main result.

%th5.1 #&#
\begin{theo} \label{samesd}
Ordered pivotal sampling and Deville's systematic sampling with the
same parameter $\pi$ induce the same sampling design.
\end{theo}

\begin{pf}%[of Theorem~\ref{samesd}]
From Lemmas~\ref{twostage} and~\ref{twostage2}, it is sufficient
to prove the result in case of ordered systematic sampling and
Deville's systematic sampling with parameter $\psi$ in the
population $U_c$. We only need to show that equations
(\ref{case1})--(\ref{case3}) hold in case of ordered pivotal
sampling. Recall that we note
\begin{eqnarray*}
U_{c,i} & = & \{u_{2i-2},\ldots,u_{2n-1}\}, \\
\psi_i & = &  (b_{i-1},\phi_{2i-1},\ldots,\phi_j,\ldots,\phi
_{2n-1} )',
\end{eqnarray*}
and that $S_{op,i}$ denotes a random sample selected in $U_{c,i}$ by
means of ordered pivotal sampling with parameter $\psi_i$ (see
Section~\ref{secops}).

We first consider equation (\ref{case1}). From (\ref{equwj1}), we
obtain:
\begin{eqnarray*}
& & pr ( X_{i+1}=u_{2i} | X_1,\ldots,X_{i-1},X_i=u_{2i-2}
) \\
&& \quad  =  pr ( X_{i+1}=u_{2i} | X_1,\ldots,X_{i-1},X_i=u_{2i-2},
J_{i-1} \in\{ X_1,\ldots,X_i \}  ),
\end{eqnarray*}
which is equivalent to $pr ( u_{2i} \in S_{op,i}, u_{2i-1}
\notin S_{op,i} | u_{2i-2} \in S_{op,i}  )$, so that the
result follows from equation (\ref{opsi1}).

Similarly, we obtain
\begin{eqnarray*}
& & pr ( X_{i+1}=u_{2i+1} | X_1,\ldots,X_{i-1},X_i=u_{2i-2}
 ) \\
&& \quad  =  pr ( X_{i+1}=u_{2i+1} | X_1,\ldots,X_{i-1},X_i=u_{2i-2},
J_{i-1} \in\{ X_1,\ldots,X_i \}  ) \\
&& \quad  \equiv pr ( u_{2i+1} \in S_{op,i},u_{2i} \notin S_{op,i},
u_{2i-1} \notin S_{op,i} | u_{2i-2} \in S_{op,i}  ) \\
&& \quad  =  %\frac{ b_{i-1} (1-a_i-b_i)  \{(1-a_i)(1-b_i) \}^{-1}
\frac{(1-a_i-b_i)(1-b_i-a_{i+1})}{(1-a_i)(1-b_i)},
\end{eqnarray*}
where the last line follows from (\ref{opsi2}), and
\begin{eqnarray*}
& & pr ( X_{i+1}=u_{2i+2} | X_1,\ldots,X_{i-1},X_i=u_{2i-2}
 ) \\
&& \quad  =  pr ( X_{i+1}=u_{2i+2} | X_1,\ldots,X_{i-1},X_i=u_{2i-2},
J_{i-1} \in\{ X_1,\ldots,X_i \}  ) \\
&& \quad  \equiv pr ( u_{2i+2} \in S_{op,i},u_{2i+1} \notin
S_{op,i},u_{2i} \notin S_{op,i}, u_{2i-1} \notin S_{op,i} | u_{2i-2}
\in S_{op,i}  ) \\
&& \quad  =  %\frac{ b_{i-1} a_{i+1} (1-a_i-b_i)  \{ (1-a_i) (1-b_i)
%(1-b_{i+1})  \}^{-1} (1-b_{i+1})}{b_{i-1}} =
\frac{(1-a_i-b_i)a_{i+1}}{(1-a_i)(1-b_i)},
\end{eqnarray*}
where the last line follows from (\ref{opsi3}). This proves
equation (\ref{case1}). The proof for equation (\ref{case3}) is
similar, and is thus omitted.

We now turn to equation (\ref{case2}). We introduce some further
notation. Let
\begin{eqnarray*}
U_{c,i+1} & = & \{u_{2i},\ldots,u_{2n-1}\}, \\
\psi_{i+1} & = &  (b_{i},\phi_{2i+1},\ldots,\phi_j,\ldots,\phi
_{2n-1} )',
\end{eqnarray*}
and let $S_{op,i+1}$ be a random sample selected in $U_{c,i+1}$ by
means of ordered pivotal sampling with parameter $\psi_{i+1}$. We
have
%
%e5.1 ###
%e5.1 #&#
\begin{eqnarray}\label{eq1}
& & pr ( X_{i+1}=u_{2i} | X_1,\ldots
,X_{i-1},X_i=u_{2i-1},W_i=u_{2i-1}  ) \nonumber\\
&& \quad  =  pr ( X_{i+1}=u_{2i} | X_1,\ldots
,X_{i-1},X_i=u_{2i-1},J_i=u_{2i}  ) \\
&& \quad  \equiv pr ( u_{2i} \in S_{op,i+1}  ) = b_i,\nonumber
\end{eqnarray}
where the second line in (\ref{eq1}) comes from
\[
\{ X_i=u_{2i-1},W_i=u_{2i-1} \} \Leftrightarrow\{
X_i=u_{2i-1},J_i=u_{2i} \}.
\]
Also,
%
%e5.2 ###
%e5.2 #&#
\begin{eqnarray}
pr ( X_{i+1}=u_{2i} | X_1,\ldots
,X_{i-1},X_i=u_{2i-1},J_i=u_{2i-1}  ) & = & 1, \label{eq2}
\end{eqnarray}
since
\begin{eqnarray*}
\{ X_i=u_{2i-1},J_i=u_{2i-1} \} & \Rightarrow& \{ X_i=u_{2i-1},W_i=u_{2i}
\}  \Rightarrow \{ X_{i+1}=u_{2i}\} .
\end{eqnarray*}

Further,
%
%e5.3 ###
%e5.3 #&#
\begin{eqnarray}\label{eq3}
& & pr ( W_{i}=u_{2i-1} | X_1,\ldots,X_{i-1},X_i=u_{2i-1}
) \nonumber\\
&& \quad  =  pr ( W_{i}=u_{2i-1} | X_i=u_{2i-1}  ) \\
&& \quad  =  pr ( X_{i}=u_{2i-1} | W_i=u_{2i-1}  ) \frac{pr (
W_{i}=u_{2i-1}  )}{pr ( X_{i}=u_{2i-1}  )}  \nonumber\\
&& \quad  =  1 \times\frac{(1-a_i-b_{i-1})(1-a_i-b_i)  \{ (1-a_i)
(1-b_i)  \}^{-1}}{1-a_i-b_{i-1}} \nonumber\\
&& \quad  =  \frac{1-a_i-b_i}{(1-a_i)(1-b_i)}, \nonumber
\end{eqnarray}
the fourth line in (\ref{eq3}) being a consequence of Lemma
\ref{probu2i-1}. The same reasoning leads to
%
%e5.4 ###
%e5.4 #&#
\begin{eqnarray}\label{eq4}
& & pr ( J_{i}=u_{2i-1} | X_1,\ldots,X_{i-1},X_i=u_{2i-1}
) \nonumber\\
&& \quad  =  pr ( J_{i}=u_{2i-1} | X_i=u_{2i-1}  ) \nonumber\\
&& \quad  =  pr ( X_{i}=u_{2i-1} | J_i=u_{2i-1}  ) \frac{pr (
J_{i}=u_{2i-1}  )}{pr ( X_{i}=u_{2i-1}  )}  \\
&& \quad  =  b_i \times\frac{a_i(1-a_i-b_{i-1})  \{ (1-a_i) (1-b_i)
 \}^{-1}}{1-a_i-b_{i-1}} \nonumber\\
&& \quad  =  \frac{a_i b_i}{(1-a_i)(1-b_i)}. \nonumber
\end{eqnarray}

From equations (\ref{eq1})--(\ref{eq4}), we obtain that
\begin{eqnarray*}
& & pr ( X_{i+1}=u_{2i} | X_1,\ldots,X_{i-1},X_i=u_{2i-1}
 ) \\
&& \quad  =  b_i \times\frac{1-a_i-b_i}{(1-a_i)(1-b_i)} + 1 \times\frac{a_i
b_i}{(1-a_i)(1-b_i)} \\
&& \quad  =  \frac{b_i}{1-a_i}.
\end{eqnarray*}
Similar computations lead to
\begin{eqnarray*}
pr ( X_{i+1}=u_{2i+1} | X_1,\ldots,X_{i-1},X_i=u_{2i-1}  )
& = & \frac{(1-b_i-a_{i+1})(1-a_i-b_i)}{(1-a_i)(1-b_i)}
\end{eqnarray*}
and
\begin{eqnarray*}
pr ( X_{i+1}=u_{2i+2} | X_1,\ldots,X_{i-1},X_i=u_{2i-1}  )
& = &
\frac{a_{i+1}(1-a_i-b_i)}{(1-a_i)(1-b_i)},
\end{eqnarray*}
which proves (\ref{case2}).
\end{pf}

Theorem~\ref{samesd} implies that ordered pivotal sampling shares the
same second-order inclusion probabilities as Deville's systematic
sampling. The computation of these probabilities is developed in  Deville~\cite{dev88}, and is reminded below.

%th5.2 #&#
\begin{theo}[(Deville~\cite{dev88})] \label{pikl}
Let $k$ and $l$ be two distinct units in $U$. If $k$ and $l$ are two
non-cross-border units that belong to the same microstratum
$U_i$, then
\[
\pi_{kl}=0,
\]
if $k$ and $l$ are two non-cross-border units that belong to distinct
microstrata $U_i$ and $U_j$, respectively, where $i<j$, then
\[
\pi_{kl}=\pi_k \pi_l  \{ 1-c(i,j) \},
\]
if $k=k_{i-1}$ and $l$ is a non-cross-border unit that belongs to the
microstratum $U_j$ where $i \leq j$, then
\[
\pi_{kl}=\pi_k \pi_l  [ 1-b_{i-1} (1-\pi_k)  \{\pi_k
(1-b_{i-1})  \}^{-1} c(i,j)  ],
\]
if $l=k_{j-1}$ and $k$ is a non-cross-border unit that belongs to the
microstratum $U_i$ where $i < j$, then
\[
\pi_{kl}=\pi_k \pi_l  \{ 1-(1-\pi_l)(1-b_{j-1}) (\pi_l
b_{j-1})^{-1} c(i,j)  \},
\]
if $k=p_{i-1}$ and $l=p_{j-1}$, where $i < j$, then
\[
\pi_{kl}=\pi_k \pi_l  [ 1- b_{i-1} (1-b_{j-1}) (1-\pi_k) (1-\pi_l)
 \{\pi_k \pi_l b_{j-1} (1-b_{i-1})  \}^{-1} c(i,j)  ],
\]
where $c(i,j)=\prod_{l=i}^{j-1} c_l$, $c_l=a_l b_l  \{
(1-a_l)(1-b_l) \}^{-1}$ and with $c(i,i)=1$.
\end{theo}

As noticed by Deville~\cite{dev88}, it follows from Theorem~\ref{pikl} that
many of the second-order inclusion probabilities are zero. As a result,
no unbiased variance estimator may be found for the Horvitz--Thompson
estimator. The search for variance estimators under reasonable model
assumptions for the variable of interest $y$ is a matter for further research.

%s6 ###
%se6 #&#
\section{Interest of ordered pivotal sampling} \label{intops}

This is clear from Theorems~\ref{samesd} and~\ref{pikl} that ordered
pivotal sampling induces a sampling design with a rather small entropy,
since the second-order inclusion probabilities heavily depend on the
order of the units in the population. If the maximization of entropy is
a major concern, \textit{randomized pivotal sampling}, where the list of
the units in the population is randomly ordered before applying the
pivotal method, should certainly be preferred. The main interest of
ordered pivotal sampling lies in the gain of precision obtained from a
stratification effect, if the ranking of the units in the population is
well correlated to the variable of interest. In this sense, ordered
pivotal sampling is similar in spirit to classical, ordered systematic
sampling. However, systematic sampling can be particularly inefficient
if the ordering is unappropriate, with regard to the variable of
interest. Ordered pivotal sampling introduces more randomization in the
sampling process, and should be more robust, in some sense, than
systematic sampling. In the sequel, ordered pivotal sampling is
compared to other sampling designs with respect to various criteria.

To fix ideas, we consider the case of (i) equal inclusion probabilities
$\pi_k=n/N$, such that (ii) the population size $N$ is an integer
multiple of the sample size $n$, and we note $N=n p$. In this case, the
microstrata $U_i, i=1,\ldots,n,$ are non overlapping with the same size
$N_i=p$. We have
%
%e6.1 ###
%e6.1 #&#
\begin{equation} \label{microstrateqprob}
U_i=\{(i-1)p+1,\ldots,(i-1)p+p\},
\end{equation}
and ordered pivotal sampling amounts to stratified simple random
sampling of size $n_i=1$ inside each microstratum $U_i$. Also, it is
well known that under the same assumptions (i) and (ii), systematic
sampling amounts to simple random sampling of size $m=1$ in the
population $G_c=\{g_1,\ldots,g_p\}$ of $M=p$ clusters, where each cluster
%
%e6.2 ###
%e6.2 #&#
\begin{equation} \label{clustsys}
g_j=\{j,j+p,\ldots,j+(n-1)p\}
\end{equation}
contains $M_j=n$ units. Let $y$ denote some variable of interest, and let
%
%e6.3 ###
%e6.3 #&#
\begin{equation} \label{horthoest}
\hat{t}_{y\pi}=\sum_{k \in S} \frac{y_k}{\pi_k}
\end{equation}
denote the Horvitz--Thompson (HT) estimator of the total $t_y=\sum_{k
\in U} y_k$.

Under conditions (i) and (ii), ordered systematic sampling and ordered
pivotal sampling may be seen as particular cases of Markov chain
designs (Breidt~\cite{bre95}). Let $M$ be a doubly stochastic transition
probability matrix, of size $p$. In a Markov chain design with matrix
of transition $M$, a sample $s=\{R_1,p+R_2,\ldots,(n-1)p+R_n\}$ is
selected, where
$R_1,\ldots,R_n$ is the Markov chain associated to $M$, with $R_1$
being uniformly distributed on $\{1,\ldots,p\}$. Let $I(p)$ denote the
identity matrix of size $p$, and $J(p)$ denote the square matrix of
size $p$ with all elements equal to $1$. The use of the matrix of transition
\[
M_{\rho}=\rho\frac{J(p)}{p} + (1-\rho) I(p),
\]
with $\rho\in[0,1]$ defines the category of compromise Markov chain
designs (Breidt~\cite{bre95}). The choice $\rho=0$ leads to ordered systematic
sampling, while the choice $\rho=1$ leads to ordered pivotal sampling.

%s6.1 ###
%su6.1 #&#
\subsection{Entropies of sampling designs}

As a measure of randomness of a sampling design $q(\cdot)$, we use the
entropy $H(q)$ defined as
%
%e6.4 ###
%e6.4 #&#
\begin{equation} \label{entropy}
H(q)=-\sum_{s \subset U} q(s) \log q(s),
\end{equation}
with $0 \log0 = 0$ by convention. We have
\begin{eqnarray*}
H(srs) & = & \log N! - \log n! - \log(N-n)!  \\
& = & \sum_{k=0}^{n-1} \log \biggl( \frac{N-k}{n-k}  \biggr)
\end{eqnarray*}
for simple random sampling, and
\begin{eqnarray*}
H(sys) & = & \log \biggl( \frac{N}{n}  \biggr)
\end{eqnarray*}
for ordered systematic sampling, see for example Till\'e and Haziza~\cite{tilhaz10}.
Some straightforward algebra leads to
\begin{eqnarray*}
H(ops) & = & n \log \biggl( \frac{N}{n}  \biggr)
\end{eqnarray*}
for ordered pivotal sampling. As a measure of comparison of entropy for
two sampling designs $q(\cdot)$ and $r(\cdot)$, we may use the
Kullback--Leibler divergence
\[
D(q \| r)=\sum_{s \subset U} q(s) \log\frac{q(s)}{r(s)}
\]
if the two sampling designs are such that $r(s)=0 \Rightarrow q(s)=0$.
We obtain
\begin{eqnarray*}
D(sys \| srs) & = & \sum_{k=1}^{n-1} \log \biggl( \frac{N-k}{n-k}
\biggr), \\
D(ops \| srs) & = & \sum_{k=0}^{n-1} \log\biggl ( \frac{1-k/N}{1-k/n}
 \biggr),\\
D(sys \| ops) & = & (n-1) \log \biggl( \frac{N}{n}  \biggr).
\end{eqnarray*}
Both simple random sampling and ordered pivotal sampling clearly have
much larger entropy than ordered systematic sampling.

%s6.2 ###
%su6.2 #&#
\subsection{Maximum design-effect for sampling designs}

This is a standard fact that the variance of the HT-estimator under
without-replacement simple random sampling is given by
%
%e6.5 ###
%e6.5 #&#
\begin{equation} \label{varsrs}
V_{srs}  ( \hat{t}_{y\pi}  ) = N^2 \frac{1-f}{n} S_y^2,
\end{equation}
where $f=n/N$, $S_{y}^2=\frac{1}{N-1} \sum_{k \in U}  (y_k - \mu
_{y}  )^2$ and $\mu_{y}=\frac{1}{N} \sum_{k \in U} y_k$. On the
other hand, the variance of the HT-estimator under ordered pivotal
sampling and assumptions (i) and (ii) may be written as
%
%e6.6 ###
%e6.6 #&#
\begin{equation} \label{varops}
V_{ops}  ( \hat{t}_{y\pi}  ) = N^2 \frac{1-f}{n} \frac{1}{n}
\sum_{i=1}^n S_{yi}^2,
\end{equation}
where $S_{yi}^2=\frac{1}{N_i-1} \sum_{k \in U_i}  (y_k - \mu_{yi}
 )^2$ and $\mu_{yi}=\frac{1}{N_i} \sum_{k \in U_i} y_k$. Finally,
the variance of the HT-estimator under systematic sampling is then
given by
%
%e6.7 ###
%e6.7 #&#
\begin{equation} \label{varsys}
V_{sys}  ( \hat{t}_{y\pi}  ) = N^2 \frac{1-f}{n} \frac{1}{n} S_Y^2,
\end{equation}
where
\begin{eqnarray*}
S_Y^2 & = & \frac{1}{M-1} \sum_{j=1}^p \biggl ( t_{yj} -
\frac{t_y}{M}  \biggr)^2 \\
& = & \frac{n^2}{p-1} \sum_{j=1}^p  ( m_{yj} - \mu_y
 )^2,
\end{eqnarray*}
with $t_{yj}=\sum_{k \in G_j} y_k$ and $m_{yj}=t_{yj}/n$. \\

As a measure of risk of a strategy combining a sampling design $q(\cdot
)$ and HT-estimation, we may use the maximum design-effect
%
%e6.8 ###
%e6.8 #&#
\begin{equation} \label{deffmax}
\mathit{DMAX}(q) = \max_{y \in\mathcal{C}} \frac{V_{q}  ( \hat{t}_{y\pi}
 )}{V_{srs}  ( \hat{t}_{y\pi}  )},
\end{equation}
where $\mathcal{C}$ denotes the set of non-constant variables of
interest (that is, containing all variables $y$ such that $S_y^2 \neq0$).

%th6.1 #&#
\begin{theo} \label{dmaxopssys}
Assume that conditions \textup{(i)} and \textup{(ii)} are satisfied. Then we have for
ordered pivotal sampling
%
%e6.9 ###
%e6.9 #&#
\begin{equation} \label{dmaxops}
\mathit{DMAX}(ops)=\frac{N-1}{N-n}
\end{equation}
and for ordered systematic sampling
%
%e6.10 ###
%e6.10 #&#
\begin{equation} \label{dmaxsys}
\mathit{DMAX}(sys)=n \frac{N-1}{N-n}.
\end{equation}
\end{theo}

\begin{pf}
For any variable $y$, it follows from a standard analysis of variance that
\[
S_y^2 = \sum_{i=1}^n \frac{p-1}{N-1} S_{yi}^2 + \sum_{i=1}^n \frac{p}{N-1}
 ( \mu_{yi} - \mu_y  )^2,
\]
so that
\[
\sum_{i=1}^n S_{yi}^2 \leq\frac{N-1}{p-1} S_y^2
\]
and the equality occurs if all the stratum means $\mu_{yi}$ are equal.
A joint application of (\ref{varsrs}) and (\ref{varops}) leads to
\[
\frac{V_{ops}  ( \hat{t}_{y\pi}  )}{V_{srs}  (
\hat{t}_{y\pi}  )} \leq\fraca{\frac{N-1}{n(p-1)} S_y^2}{S_y^2} =
\frac{N-1}{N-n},
\]
which gives (\ref{dmaxops}). The use of an alternative analysis
of variance leads to
\[
S_y^2 = \sum_{j=1}^p \frac{n-1}{N-1} \sigma_{yj}^2 + \sum_{j=1}^p \frac{n}{N-1}
 ( m_{yj} - \mu_y  )^2,
\]
where $\sigma_{yj}^2=\frac{1}{n-1} \sum_{k \in G_j}  (y_k -
m_{yj}  )^2$. This leads to
\[
S_Y^2 \leq\frac{n^2}{p-1} \frac{N-1}{n} S_y^2,
\]
and the equality occurs if the variable $y$ is constant inside any
cluster $g_j$. By a joint application of (\ref{varsrs}) and (\ref{varsys}),
we have
\[
\frac{V_{sys}  ( \hat{t}_{y\pi}  )}{V_{srs}  (
\hat{t}_{y\pi}  )} \leq\fraca{\frac{n^2}{p-1} \frac{N-1}{n^2}
S_y^2}{S_y^2} =
n \frac{N-1}{N-n},
\]
which gives (\ref{dmaxsys}).
\end{pf}

If the sample size $n$ remains small to moderate, equation (\ref
{dmaxops}) implies that $\mathit{DMAX}$ tends to 1 in case of ordered pivotal
sampling, if $N$ is sufficiently large. Even in the worst cases,
ordered pivotal sampling will thus be competitive to simple random
sampling. On the other hand, equation (\ref{dmaxsys}) implies that a
strategy involving systematic sampling may be considerably more risky
in some situations.

%s6.3 ###
%su6.3 #&#
\subsection{Equality of treatment of variables for sampling designs}

As pointed out by a referee, the $\mathit{DMAX}$ criterion considered previously
is very stringent since referring to the worst possible variable for
sampling designs. In case of ordered systematic sampling, this would be
a cyclical variable whose period is equal to $p=N/n$; such a situation
is usually unlikely to occur, except in particular situations.

An alternative criterion studied by Deville~\cite{dev87} and Qualit\'e~\cite{qua10}
considers the equality of treatment of variables. For any sampling
design $q(\cdot)$ with first-order inclusion probabilities $\pi_k$ and
second-order inclusion probabilities $\pi_{kl}(q)$ with $k,l \in U$, let
\[
\Delta(q)= [\Delta_{kl}(q) ]_{k,l \in U}
\]
be its design variance--covariance matrix, with $\Delta_{kl}(q)=\pi
_{kl}(q)-\pi_k \pi_l$. Let $0 \leq\lambda_1(q) \leq\cdots \leq\lambda
_N(q)$ denote the eigenvalues of $\Delta(q)$. For any variable $y$ that
lies on the unit sphere of the Euclidean norm (that is, such that $\sum
_{k \in U} y_k^2=1$), we have
\[
\lambda_1(q) \leq V_p  ( \hat{t}_{y\pi}  ) \leq\lambda_N(q).
\]
Roughly speaking, the extreme eigenvalues give the extreme possible
values for the variance. Qualit\'e~\cite{qua10}, page 50, proposed to measure the
equality of treatment for variables in terms of minimization of the
dispersion of the eigenvalues, denoted as
\[
\delta(q)=N^{-1} \sum_{k=1}^N  \{\lambda_k(q) - \bar{\lambda}(q)
 \}^2
\]
with $\bar{\lambda}(q)=N^{-1} \sum_{k=1}^N \lambda_k(q)$. Note that
$\sum_{k=1}^N \lambda_k(q)=Tr(\Delta(q))$, where $Tr(\cdot)$ denotes
the trace. For any sampling design $q(\cdot)$ with equal probabilities
$\pi_k=n/N$, this leads to
\begin{eqnarray*}
\bar{\lambda}(q) & = & N^{-1} \sum_{k \in U} \pi_k(1-\pi_k) \\
& = & \frac{p-1}{p^2},
\end{eqnarray*}
which will be simply denoted as $\bar{\lambda}$ in the sequel.

The ranking of the evaluated sampling designs with respect to this
criterion is established in Theorem~\ref{equtreat}. Clearly, ordered
pivotal sampling tends to treat the variables more equally than ordered
systematic sampling. To demonstrate Theorem~\ref{equtreat}, we need
the following lemma.

%le6.1 #&#
\begin{lemma} \label{delta2vp}
Assume that the sampling design $q(\cdot)$ is performed with equal
probabilities $\pi_k=n/N$, and
that conditions \textup{(i)} and \textup{(ii)} are satisfied. Assume that
$\Delta(q)$ has only two eigenvalues $0$ and $\lambda_{+}(q)>0$,
with multiplicities $N_0(q)$ and $N-N_0(q)$,
respectively. Then:
%
%e6.11 ###
%e6.11 #&#
\begin{equation} \label{delta2vpform}
\delta(q)=\frac{N_0(q)}{N-N_0(q)} \bar{\lambda}^2.
\end{equation}
\end{lemma}

\begin{pf}
Since $q(\cdot)$ has only one strictly
positive eigenvalue $\lambda_{+}(q)$ with multiplicity
$N-N_0(q)$, we have $\lambda_{+}(q)=\frac{N}{N-N_0(q)}
\bar{\lambda}$. Then
\begin{eqnarray*}
N \delta(q) & = & \sum_{k=1}^N  \{\lambda_k(q) - \bar{\lambda}
 \}^2
\\
& = & N_0(q) \bar{\lambda}^2 + \{N-N_0(q)\} \biggl \{\frac{N}{N-N_0(q)}
\bar{\lambda} - \bar{\lambda}  \biggr\}^2 \\
& = & N \frac{N_0(q)}{N-N_0(q)} \bar{\lambda}^2.
\end{eqnarray*}
\upqed
\end{pf}

%th6.2 #&#
\begin{theo} \label{equtreat}
Assume that conditions \textup{(i)} and \textup{(ii)} are satisfied. Then we have
%
%e6.12 ###
%e6.12 #&#
\begin{equation} \label{equtreatform}
\delta(srs) \leq\delta(ops) \leq\delta(sys).
\end{equation}
\end{theo}

\begin{pf}
It may be easily shown (see, e.g., Deville~\cite{dev87}, page 120) that in
case of simple random sampling, $\Delta(srs)$ has
only two eigenvalues $0$ and $\lambda_{+}(srs)=\frac{N \bar{\lambda}}{N-1}$,
with \mbox{$N_0(srs)=1$}.

Let $I(n)$ denotes the identity matrix of size $n$, and $J(n)$ denotes
the square matrix of size
$n$ with all elements equal to $1$. After some algebra, we obtain from
the definition of ordered pivotal sampling that
\[
\Delta(ops)= \pmatrix{
\Delta_1 & 0 &\cdots & 0 \cr
0 & \Delta_1 &\cdots & \vdots\cr
\vdots&\ddots & \ddots& 0 \cr
0 &\cdots & 0 & \Delta_1
}=I(n) \otimes\Delta_1,
\]
where $\Delta_1=p^{-1} \{ I(p)-p^{-1} J(p) \}$ and $\otimes$ denotes
the Kronecker product. It follows that the $N$ eigenvalues of $\Delta
(ops)$ are given by the products of the eigenvalues of $I(n)$ and
$\Delta_1$ (see for example Theorem 1 in Magnus and Neudecker~\cite{magneu88}, page 28). The
eigenvalues of $\Delta_1$ are $p^{-1}$ and $0$ with multiplicities
$p-1$ and $1$, respectively. Consequently, $\Delta(ops)$ has only two
eigenvalues $0$ and $\lambda_{+}(ops)=p^{-1}$, with $N_0(ops)=n$.

Similarly, we obtain from the definition of ordered systematic
sampling that
\[
\Delta(sys)= \pmatrix{
\Delta_1 &\cdots &\cdots & \Delta_1 \cr
\vdots& & & \vdots\cr
\vdots& & & \vdots\cr
\Delta_1 &\cdots &\cdots & \Delta_1 }
 =J(n) \otimes\Delta_1.
\]
Since the eigenvalues of $J(n)$ are $0$ and $n$ with multiplicities
$n-1$ and $1$, respectively,
$\Delta(sys)$ has only two eigenvalues $0$ and $\lambda_{+}(sys)=n p^{-1}$,
with $N_0(sys)=N-p+1$.

Equation (\ref{delta2vpform}) implies that $\delta(q)$
increases as $N_0(q)$ increases. Clearly, $N_0(srs) \leq N_0(ops)$, and
from the identity $(N-p+1)-n=(p-1)(n-1) \geq0$, we obtain that
$N_0(ops) \leq
N_0(sys)$ so that the result follows.
\end{pf}

%s6.4 ###
%su6.4 #&#
\subsection{Some numerical results on a small population}

To investigate on the properties of considered sampling algorithms, we
considered a small example. We first generated a finite population of
size $N=12$, containing three variables of interest, $y_1$, $y_2$ and
$y_3$. Table~\ref{yvalues} shows the values for the three variables of
interest. The variable $y_1$ is highly correlated to the order of the
units in the population, on the contrary to variable $y_2$. The
variable $y_3$ exhibits a particularly unfavorable case for systematic sampling.

%t2 ###
%ta2 #&#
\begin{table*}
\tabcolsep=0pt
\caption{Values of three variables of interest in the generated
population} \label{yvalues}
\begin{tabular*}{\textwidth}{@{\extracolsep{\fill}}lllllllllllll@{}}
\hline
Unit & $1$ & $2$ & $3$ & $4$ & $5$ & $6$ & $7$ & $8$ & $9$ & $10$ &
$11$ & $12$ \\
\hline
$y_1$ & 10 & 10 & 10 & 15 & 45 & 45 & 50 & 50 & 60 & 60 & 60 & 65 \\
$y_2$ & 15 & 45 & 10 & 60 & 60 & 50 & 45 & 65 & 10 & 50 & 10 & 60 \\
$y_3$ & 10 & 45 & 60 & 15 & 50 & 65 & 10 & 50 & 60 & 10 & 45 & 60 \\
\hline
\end{tabular*}
\end{table*}

We considered equal probability sampling of size $n=2$ (respectively,
$n=4$) by means of six sampling designs: simple random sampling
without replacement (SRS), ordered systematic sampling (SYS),
compromise Markov chain design with $\rho=0.25$ (CMC25), $\rho=0.50$
(CMC50), $\rho=0.75$ (CMC75), and ordered pivotal sampling (OPS). As a
measure of variability of the HT-estimator $\hat{t}_{y\pi}$ for a
sampling design $q(\cdot)$, we considered the design-effect (DEFF)
given by
%
%e6.13 ###
%e6.13 #&#
\begin{equation} \label{deff}
\mathit{DEFF}=\frac{V_{q}  ( \hat{t}_{y\pi}  )}{V_{srs}  ( \hat
{t}_{y\pi}  )},
\end{equation}
where the variances are computed by means of formulas (\ref
{varsrs})--(\ref{varsys}) for SRS, OPS and SYS, and from formula (1)
in Breidt~\cite{bre95}, page 66, for compromise Markov chain designs. Table \ref
{simdeff} shows DEFF for the five strategies. As could be expected,
the CMC25, CMC50 and CMC75 give compromise results between SYS and OPS.
Also, it is clear from Table~\ref{simdeff} that both OPS and SYS lead
to a subsequent reduction of variance for variable $y_1$, with $\mathit{DEFF}$
ranging from $0.17$ to $0.50$ and OPS performing significantly better.
The OPS strategy is essentially similar to SRS for the variable $y_2$
which is poorly correlated to the order of the units in the population,
while SYS may be much worse ($\mathit{DEFF}=1.39$ for $n=2$) or much better
($\mathit{DEFF}=0.36$ for $n=4$). Finally, we obtain for the variable $y_3$ a
considerable loss for SYS, while the loss is more limited for OPS
with $\mathit{DEFF}=1.10$ for $n=2$ and $\mathit{DEFF}=1.36$ for $n=4$.

%t3 ###
%ta3 #&#
\begin{table*}
\tabcolsep=0pt
\caption{Design-effect for three variables of interest and five
strategies in the generated population} \label{simdeff}
\begin{tabular*}{\textwidth}{@{\extracolsep{4in minus 4in}}lllllll@{}}
\hline
&  \multicolumn{3}{l}{Sample size $n=2$} &  \multicolumn{3}{l@{}}{Sample size $n=4$}\\[-5pt]
&  \multicolumn{3}{l}{\hrulefill} &  \multicolumn{3}{l@{}}{\hrulefill}\\
&   $y_1$   & $y_2$   & $y_3$   & $y_1$   & $y_2$   & $y_3$\\
\hline
SYS & 0.50 & 1.39 & 2.18 & 0.27 & 0.36 & 5.44 \\
CMC25 & 0.46 & 1.31 & 1.91 & 0.24 & 0.61 & 3.94 \\
CMC50 & 0.43 & 1.24 & 1.64 & 0.21 & 0.76 & 2.81 \\
CMC75 & 0.39 & 1.17 & 1.37 & 0.19 & 0.85 & 1.97 \\
OPS & 0.35 & 1.10 & 1.10 & 0.17 & 0.95 & 1.36 \\
\hline
\end{tabular*}
\end{table*}

\section*{Acknowledgements}
The author wishes to thank two referees for useful comments and
suggestions that helped improving the quality of the paper significantly.

% imsref loaded by smiklovaite, 2011-11-08 14:00:53

\printhistory

\end{document}